\documentclass[twocolumn]{article}
\usepackage{nolta2017}
\usepackage{txfonts}

\usepackage[superscript,sort]{cite}
\usepackage{ifthen,ifpdf}
\ifpdf
	\usepackage[pdftex]{hyperref}	 \hypersetup{colorlinks=true,linkcolor=black,citecolor=black,urlcolor=blue,backref=page,bookmarks=true,breaklinks=true,plainpages=false }
	\usepackage[pdftex]{graphicx}
	\usepackage[usenames,pdftex]{color}
\else
	 \usepackage[colorlinks=true,breaklinks=true,plainpages=false]{hyperref}
	\usepackage{graphicx}
	\usepackage[usenames]{color}
\fi
\usepackage[textwidth=0.8in,textsize=footnotesize]{todonotes}

\begin{document}

\title{Deep Learning for Real Time Crime Forecasting}

\author{Bao Wang${}^\dag$*\thanks{*corresponding author}, Duo Zhang${}^\dag$, Duanhao Zhang${}^\dag$, P. Jeffrey Brantingham${}^\ddag$, and Andrea L. Bertozzi${}^\dag$}
\address{
\dag Department of Mathematics, UCLA\\
\ddag Department of Anthropology, UCLA\\
Email: wangbaonj@gmail.com, jasonzhang56886@gmail.com, adylucioforzm@live.com,\\
branting@g.ucla.edu, bertozzi@math.ucla.edu
}
\maketitle

\abstract
Accurate real time crime prediction is a fundamental issue for public safety, but remains a challenging problem for the scientific community. Crime occurrences depend on many complex factors. Compared to many predictable events, crime is sparse. At different spatio-temporal scales, crime distributions display dramatically different patterns. These distributions are of very low regularity in both space and time. In this work, we adapt the state-of-the-art deep learning spatio-temporal predictor, ST-ResNet [Zhang et al, AAAI, 2017], to collectively predict crime distribution over the Los Angeles area. Our models are two staged. First, we preprocess the raw crime data. This includes regularization in both space and time to enhance predictable signals. Second, we adapt hierarchical structures of residual convolutional units to train multi-factor crime prediction models. 
Experiments over a half year period in Los Angeles reveal highly accurate predictive power of our models.

\endabstract


\section{Introduction}
Real time crime forecasting is an important scientific and sociological problem. It is directly related to our quality of life. Recent efforts have been devoted to the mathematical modeling of crime. Short et al developed novel differential equation models for dynamics of crime hotspots \cite{Short:2008M3AS}. Considering crime as self-exciting, Mohler et al applied the classical epidemic type aftershock sequence (ETAS) model and its variants to crime modeling\cite{Mohler:2011JASA, Short:2014DCDSB}. These types of models provide a microscopic representation of the crime events with predictive power. The aforementioned models are built from historical data. There is also interesting work on crime prediction using social network data, e.g., Twitter \cite{XiaofengWang:2012ICSC}.

Deep learning, a cutting-edge technology for automatic feature identification via a deep neural network (DNN), gives state-of-the-art performance on many predictive scenarios, such as image classification, computer vision, speech recognition\cite{LeCun:2015Nature}. Moreover, recurrent neural networks (RNN), which are superior in mining long-range dependencies, enables highly accurate prediction of sequential data \cite{Hochreiter:1997NC}. Recently, deep learning paradigms have been applied to study spatio-temporal data. There are two frameworks of interest. One is the convolutional neural network (CNN) for learning spatial traits and RNN to learn temporal features \cite{Junbo:2017}. The other is to embed spatio-temporal data into an RNN framework via graphical models \cite{Jain:2016}.

There are many challenges in accurate real time spatio-temporal crime prediction. 
For example, compared to the traffic flow data studied in \cite{Junbo:2017}, crime data has much less spatio-temporal regularity, i.e., the number of events in adjacent time intervals and spatial cells differs hugely \cite{RN3516}. Furthermore, crime types are diverse. Our contribution in this work is three-fold. First, we select the appropriate spatio-temporal scale at which crimes are predictable. Second, we provide different approaches for data regularization in both the spatial and temporal dimensions. Third, we adapt deep learning architectures for predictions. Both convolutional and non-convolutional models are considered. In the first case, spatial connectivity is included. In the second case, grids are independent of each other.

We structure this paper as follows: In section \ref{Data}, we discuss data sets. In section \ref{Model}, the feature and deep learning architectures are presented. In section \ref{Results}, we present an example of spatio-temporal crime prediction over the last two weeks of 2015 for the whole Los Angeles (LA) area. Future work is discussed in section \ref{Conclusion}.

\section{Data Set and Preprocessing}
\label{Data}
\begin{figure}
\centering
\begin{tabular}{cc}
\includegraphics[width=0.45\columnwidth]{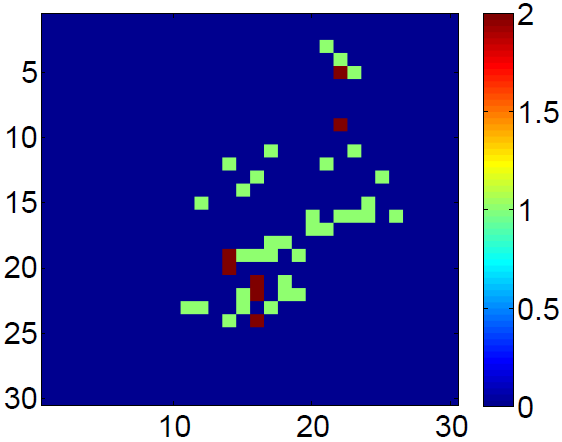}&
\includegraphics[width=0.45\columnwidth]{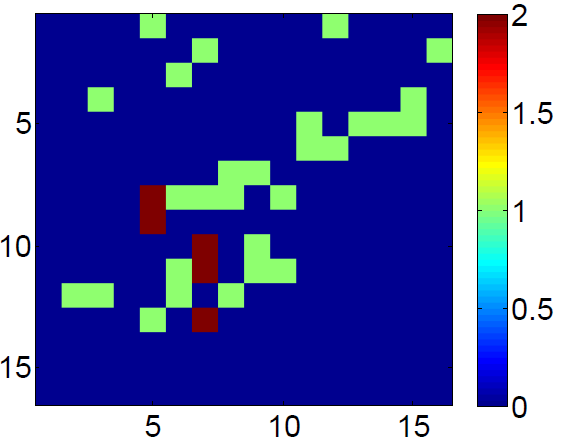}\\
(a)&(b)
\end{tabular}
\caption{Crime distribution at 1:00 p.m, Dec 20th, 2015. Chart (a) depicts crime distribution over the whole LA area; chart (b) depicts crime distribution over the selected region. The units are described in Sec. \ref{Data}.}
\label{Crime_Distribution_Demo}
\end{figure}

\paragraph{Data set description}
In this work, we consider all types of crime in LA over the last six months of 2015. In total there were 104,957 crimes. Each crime record includes crime start, end times and location. To avoid ambiguity, we regard start time of each event as the time slot.  Geographically, the latitude and longitude intervals spanned by these crimes are $[33.3427^\circ, 34.6837^\circ]$ and $[-118.8551^\circ, -117.7157^\circ]$, respectively. The spatial distribution is highly heterogeneous. A large portion of the area contains only a little crime. Therefore, we only consider crimes that happened within the region  $[33.6927^\circ, 34.3837^\circ]\times[-118.7051^\circ, -118.1157^\circ]$. In our study, we partition this selected region into a $16\times 16$ lattice.
Figure \ref{Crime_Distribution_Demo} shows the crime distribution at 1:00 p.m on Dec 20th, 2015. The left panel is the crime distribution over the whole LA area. The right panel depicts the crimes in the restricted region. The restricted region contains more than 95 percent of all crimes.

\begin{figure}
\centering
\begin{tabular}{cc}
\includegraphics[width=0.45\columnwidth]{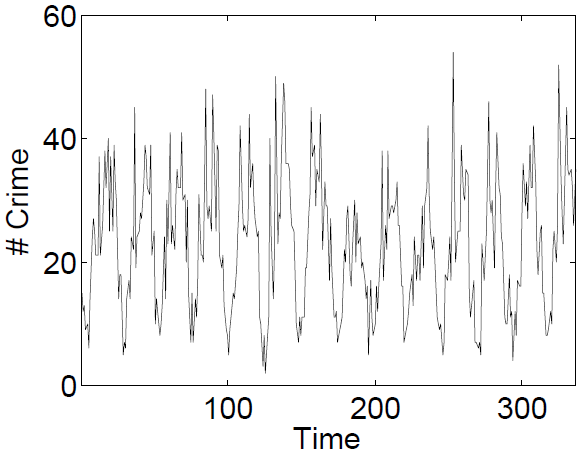}&
\includegraphics[width=0.45\columnwidth]{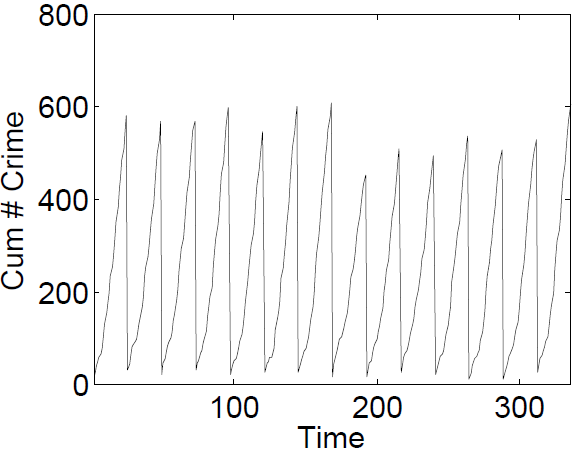}\\
(a)&(b)
\end{tabular}
\caption{Chart (a) depicts the hourly crime intensity of the last two weeks of 2015 over the whole LA area; chart (b) draws the cumulated crime intensity corresponding to (a).  Units: x-axis: time; y-axis: number of crimes.}
\label{Crime_Intensity_Demo1}
\end{figure}

\vskip -1.5cm
\paragraph{Data preprocessing}
Figure \ref{Crime_Intensity_Demo1} reflects low regularity in the temporal dimension. Nevertheless, the
hourly crime intensity indicates a strongly predictable signal; this signal is indeed periodic with a period of 24 hours.
This periodic temporal pattern exists in each grid cell. It is more spiky for the grid-wise crime intensity (see Fig.\ref{Predicted_Crime_TimeSeries}). Despite the predictability, the data is still inappropriate for deep learning which requires some regularity of the inputs, especially for grid-wise intensity. To address this, we map the original crime intensity $\{X(t)\}$ to $\{Y(t)\}$ via a periodic integral mapping:
\begin{equation}
\label{cum}
Y(t)=\int_{nT}^t X(s)ds,
\end{equation}
for $t$ within the time interval $(nT, (n+1)T]$.
Figure \ref{Crime_Intensity_Demo1} (a) shows the crime intensity over the whole LA region. 
Figure \ref{Crime_Intensity_Demo1} (b) plots the cumulative intensity of (a) after the mapping defined by Eq.(\ref{cum}).

\begin{figure}
\centering
\begin{tabular}{cc}
\includegraphics[width=0.35\columnwidth]{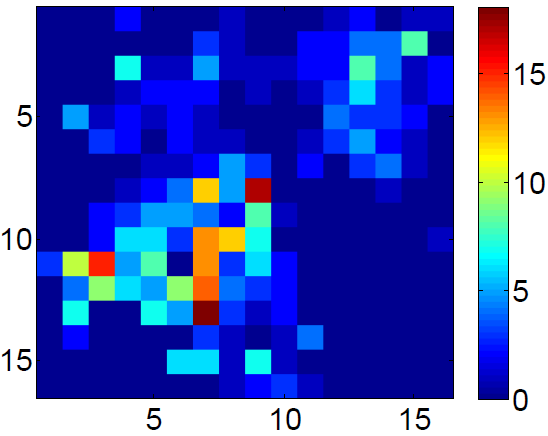}&
\includegraphics[width=0.35\columnwidth]{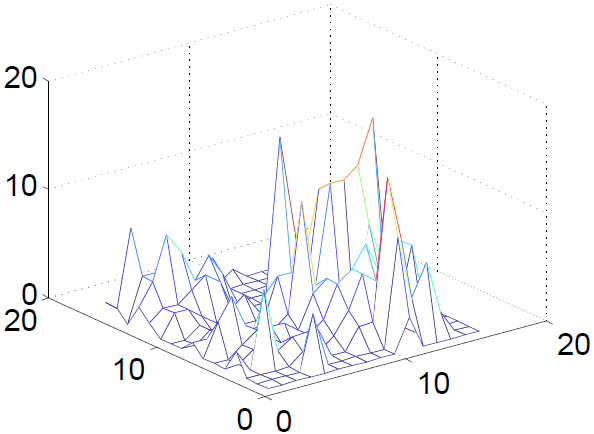}\\
(a)&(b)\\
\includegraphics[width=0.35\columnwidth]{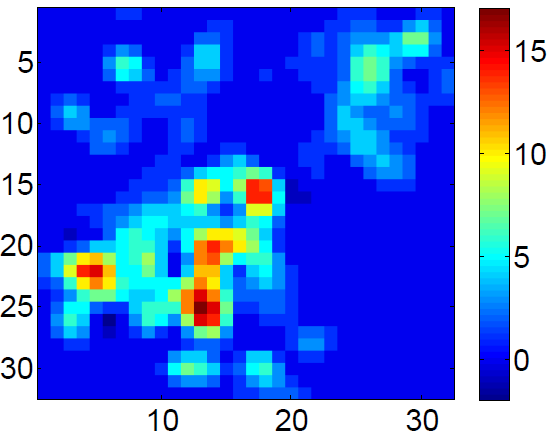}&
\includegraphics[width=0.35\columnwidth]{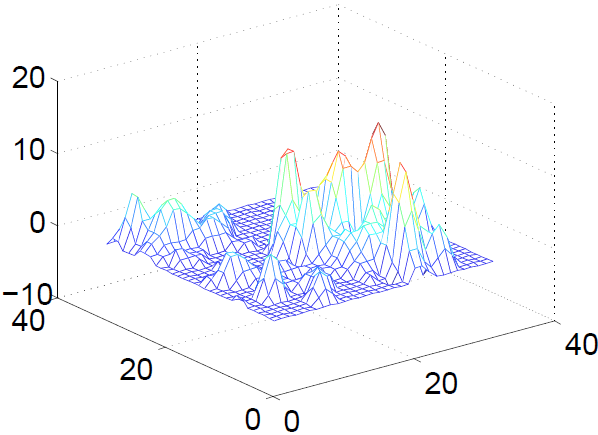}\\
(c)&(d)
\end{tabular}
\caption{Cumulated crime intensity at 11:00 p.m, Dec 31st, 2015. Chart (a) depicts crime distribution over the selected area; chart (b) provides the mesh plot of the chart (a); chart (c) depicts super resolution version of chart (a); and char (d) is mesh plot of chart (c).}
\label{Crime_cumDistribution_Demo}
\end{figure}

To address lack of spatial regularity, we can use a super resolution technique involving cubic spline interpolation. For computational efficiency, we resolve by a factor of 2 in each dimension. From Fig. \ref{Crime_cumDistribution_Demo}, we see that this cubic spline super resolution significantly improves spatial regularity. A merit of this preprocessing is that it improves the signal without losing information associated with the crime data.

\section{Models} \label{Model}
\paragraph{Problem formulation}

Given the historical filtration $\{(X_t, E_t)\}_{t=1, 2, \cdots n}$ and $\{E_{n+1}\}$, to predict $X_{n+1}$, where
$X_1, X_2, \cdots, X_{n+1}$ are the tensors representing the crime spatial distributions at times $t_1, t_2, \cdots, t_{n+1}$, $E_1, E_2, \cdots, E_{n+1}$ are external
features that affect the crimes, (e.g., holiday, time, weather), our protocol is divided into the following three steps:

\begin{itemize}
\vskip-0.5cm
\item Apply the spatio-temporal regularization to the historical data to get $\{(\hat{X}_t, E_t)\}_{t=1, 2, \cdots, n}$.
\vskip-0.5cm
\item Predict regularized cumulated crime intensity $\hat{X}^p_{n+1}$.
\vskip-0.5cm
\item Perform the spatio-temporal inverse map to get the crime intensity prediction $X^p_{n+1}$.
\vskip-0.5cm
\end{itemize}

\begin{figure}[!ht]
\small
\centering
\includegraphics[width=7cm,height=6cm]{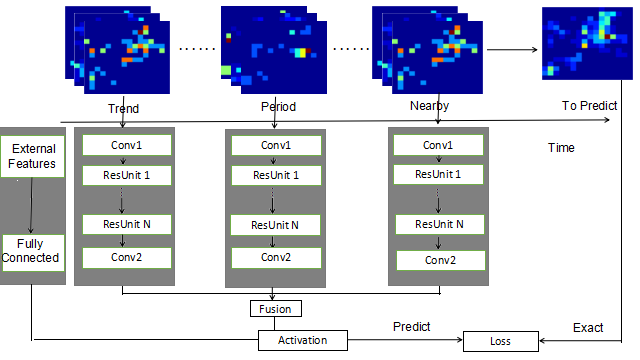}
\caption{Structure of the deep neural network models.}
\label{ConvResNet}
\end{figure}

\paragraph{Deep neural network structure}
The deep neural network architecture is depicted in Fig. \ref{ConvResNet}. We test two different structures: one is plotted in Fig. \ref{ConvResNet}; the other one is similar but without convolutional layers. For the second model, we apply an ensemble of ResNet to learn the time series on each grid, without considering the transition of crimes between different grids. The first model is more realistic. Through convolutional layers, crime dynamics can be captured. In both networks, all features are fused with the crime data itself via a parametric-matrix based fusion technique used in \cite{Junbo:2017}. The detailed description of the network structure can be found in \cite{Junbo:2017}. We implement our models using Keras \cite{Keras:2015} on top of Theano \cite{Theano:2016} software.



Our models incorporate external features, i.e., weather, holidays. Due to the periodic pattern and self exciting property of crimes \cite{Mohler:2011JASA}, we adopt trend, periodic, and nearby features. The time spacing of these features are at weekly, daily, and hourly levels, respectively.

\section{Numerical Results}\label{Results}
\begin{figure}
\centering
\begin{tabular}{cc}
\includegraphics[width=0.35\columnwidth]{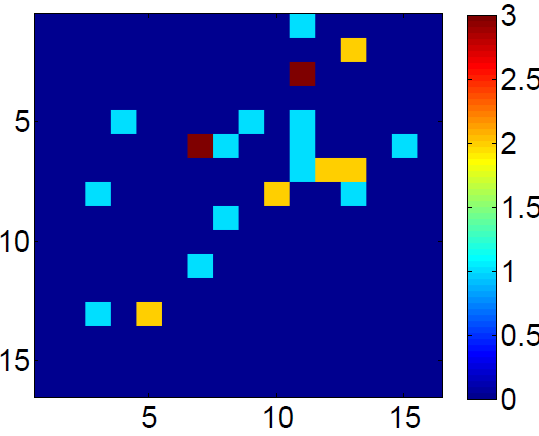}&
\includegraphics[width=0.35\columnwidth]{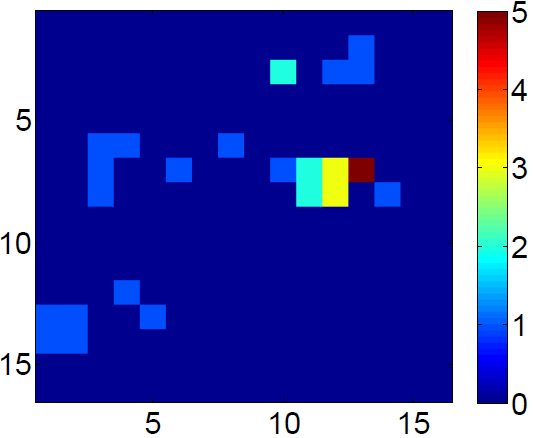}\\
(a)&(b)\\
\includegraphics[width=0.35\columnwidth]{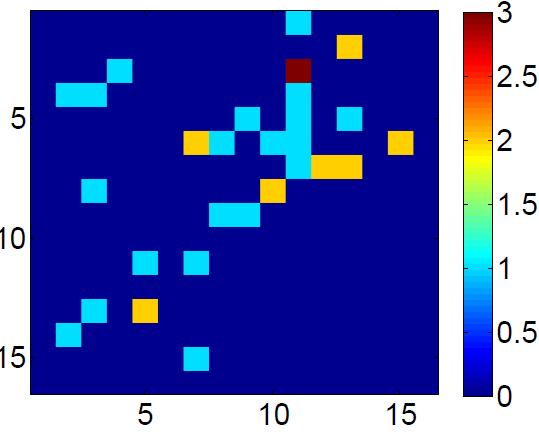}&
\includegraphics[width=0.35\columnwidth]{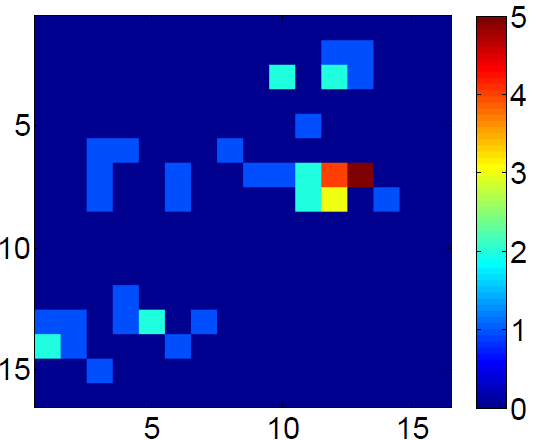}\\
(c)&(d)\\
\includegraphics[width=0.35\columnwidth]{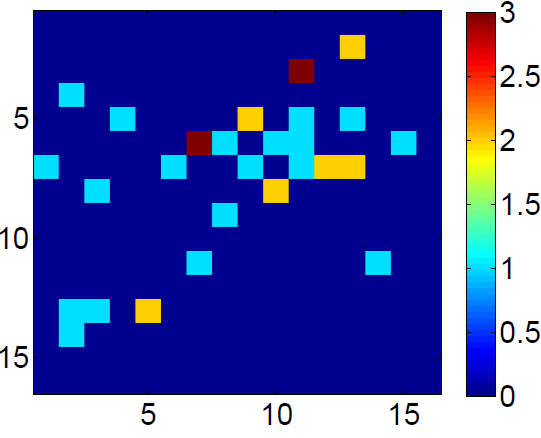}&
\includegraphics[width=0.35\columnwidth]{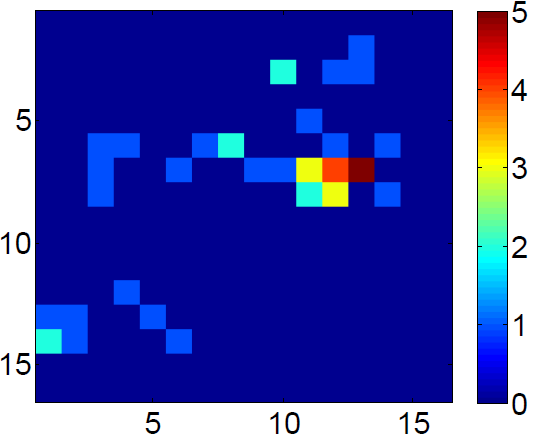}\\
(e)&(f)
\end{tabular}
\caption{Predicted vs. exact crime spatial distribution. Panels (a), (b) plot the crime spatial distribution at 1 p.m. of Dec 19, 27, 2015, respectively. Panels (c), (d) are the predicted results without convolution layers. (e), (f) are the predicted results with convolution layers.}
\label{Predicted_Crime_Distribution}
\end{figure}

In our tests, for both DNNs, the parameters are set as follows: The size of the convolution filters are fixed to $3\times 3$. The length of the trend, period, and nearby dependencies are all chosen to be 3. Six layers of residual units are used in each dependency. In the training stage, the cross validation ratio is $0.1$. The number of epochs at cross validation and fine tune periods are both 20. The learning rate is chosen to be $0.0005$. Batch normalization is used. We adopt the ADAM optimizer to optimize the loss function.
We select the last two weeks of 2015 as the test set; all the other data is used to train the DNN models.
We show an example of snapshots in time in Fig.\ref{Predicted_Crime_Distribution}.
The predictions match quite well with the ground truth. All crime hotspots are correctly predicted.

\begin{figure}
\centering
\begin{tabular}{cc}
\includegraphics[width=0.45\columnwidth]{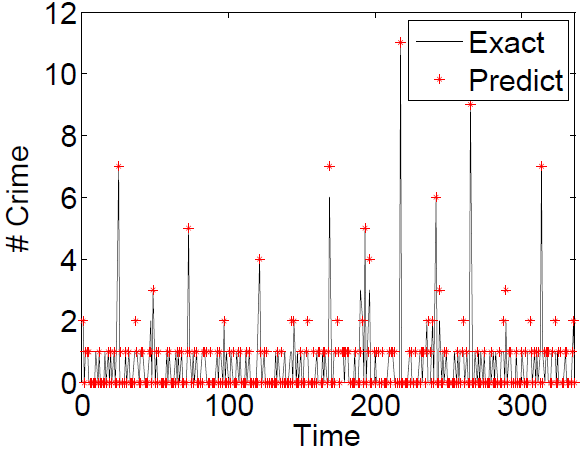}&
\includegraphics[width=0.45\columnwidth]{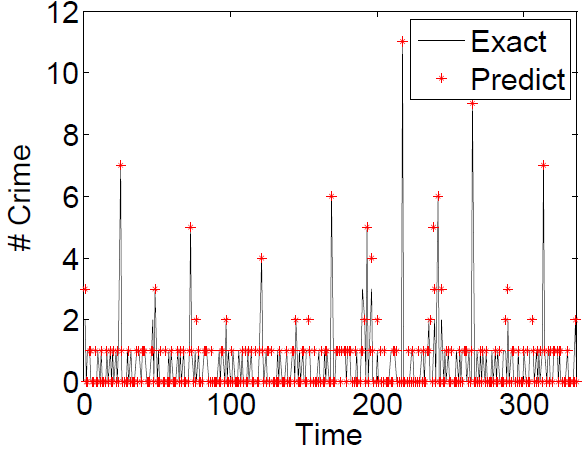}\\
(a)&(b)
\end{tabular}
\caption{Predicted vs. exact crime intensity in the area $[33.9519, 33.9951]\times [-118.3003, -118.2635]$ over the last two weeks of 2015. Charts (a) and (b) are results for models with and without convolutional layers, respectively. Units: x-axis: time; y-axis: number of crimes.}
\label{Predicted_Crime_TimeSeries}
\end{figure}

For a given grid, the crime intensity over a given time interval is also accurately predicted. As shown in Fig. \ref{Predicted_Crime_TimeSeries}, the maximum difference between the ground truth and the prediction is 2 crimes.

To quantitatively measure the performance of the models, we use the root mean squared error (RMSE):
\begin{equation}
\label{RMSE}
{\rm RMSE}=\sqrt{\frac{1}{N*T}\sum_{i, t}(I_{it}-I_{it}^p)^2}
\end{equation}
where $N$ is the total number of grids that we partition the restricted area into, $T$ is the number of time slots
considered, $I_{it}, I_{it}^p$ are the exact and predicted crime intensity in grid $i$ at time $t$, respectively. By using the same parameters as above,
the training and testing RMSEs are listed in Table \ref{RMSE}. Using the convolution technique to incorporate spatial connectivity gives better predictions.

\begin{table}[!htt]
\centering
\caption{Performance comparison between two models. Units for Train and Test columns: Number of crimes.}
\begin{tabular}{lllll}
\cline{1-5}
\footnotesize{Model} &\footnotesize{Train} &\footnotesize{Test}  &\footnotesize{$\#$parameter}s &\footnotesize{Training Time}\\
\hline
\footnotesize{Convolution}        &0.175  &0.184  &1129343 &19539s \\
\footnotesize{No Convolution}     &0.187  &0.195  &140159  &4449s \\
\hline
\end{tabular}
\label{RMSE}
\end{table}

For real applications, we focus more on the accuracy of the top N predictions. Let the top N crime hotspots be indexed as $I_N=\{i_1, i_2, \cdots, i_N\}$, where $i_1, i_2, \cdots, i_N$ are indices of the grids in the previous partition. Furthermore, denote the predicted indices of the top N crime hotspots as $I_N^p=\{j_1, j_2, \cdots, j_N\}$. We define the following accuracy of top N ranking:

\begin{equation}
\label{topN}
{\rm Acc_N}=\frac{|I_N\bigcap I_N^p|}{|I_N|},
\end{equation}
where $|*|$ is the cardinality of the set $*$. This top N recommendation system provides a practical guidance to control crimes efficiently. Table \ref{TopNACC} lists the top N accuracy averaged among all testing time slots. Using convolution provides slightly better results for all N, though the system without convolution also gives good recommendations, indicating crime hotspots are quite stationary during this test time period.
\vskip -0.5cm
\begin{table}[!htt]
\centering
\caption{Top N prediction accuracy of the above two models. Unit: $\%$.}
\begin{tabular}{llllll}
\cline{1-6}
\footnotesize{Model} &5 &10 &15 &20  &25 \\
\hline
\footnotesize{Convolution}        &80.77  &83.90 &84.35 &84.60 &84.78   \\
\footnotesize{No Convolution}
&80.48  &82.26 &82.90     &83.29 &83.95   \\
\hline
\end{tabular}
\label{TopNACC}
\end{table}
\section{Concluding Remarks} \label{Conclusion}
In this work, we apply the ST-ResNet to real time crime prediction on an hourly timescale. Due to the low regularity of the crime data in both space and time, we perform both spatial and temporal regularization of the data. More specifically, in the temporal dimension, we compute the diurnal cumulative crime per spatial region. In the spatial dimension, we use cubic spline interpolation super resolution. Compared to applying the residual network on each individual grid, our data preprocessing and CNN approach improves the prediction accuracy dramatically. The improvement is due to the fact that, in the convolutional model, the spatial information is no longer isolated allowing the crime hotspot transitions to be captured. Our predictions are extremely accurate in both space and time, which can provide reliable guidance for crime control.

There are still many things to do. In our current work, we have not classified the crime types. In reality, different types of crime may lead to different consequences. Intervention strategies should not be based on the crime intensity alone; crime type should also be addressed. Due to the super resolution regularization in space, the computational cost increases dramatically. One possible way to improve computational efficiency is to use an adaptive super resolution to train each layer of the DNN. Alternatively, we can improve the loss function of the deep learning to make the ResNet applicable to low regularity data. Another deficiency of the presented models comes from the ad hoc grid partitioning of the LA area. A better approach would be to use a graph to model the intrinsic geographic connectivity.

\section*{Acknowledgments}
We thank Dr. Penghang Yin for fruitful discussions.
This work was supported by ONR grant N00014-16-1-2119, ARO MURI grant W911NF-11-1-0332, NSF grants DMS-1118971 and DMS-1417674.
The authors thank the Los Angeles Police Department for providing the crime data for this paper.

\end{document}